%% file: FrContPal.tex
\magnification=1200

\def\Q{{\bf {Q}}}

\def\Z{{\bf Z}}     
\def\R{{\bf R}}

\input typpo

\catcode`@=11
\def\bibliographym@rk{\bgroup}
%
%
\outer\def\bye{ 	\par\vfill\supereject\end}

\def\uU{{\overline U}}

\def\uW{{\overline W}}

\def\house#1{\setbox1=\hbox{$\,#1\,$}%
\dimen1=\ht1 \advance\dimen1 by 2pt \dimen2=\dp1 \advance\dimen2 by 2pt
\setbox1=\hbox{\vrule height\dimen1 depth\dimen2\box1\vrule}%
\setbox1=\vbox{\hrule\box1}%
\advance\dimen1 by .4pt \ht1=\dimen1
\advance\dimen2 by .4pt \dp1=\dimen2 \box1\relax}

  \def\eps{{\varepsilon}}

  \def\noi{\noindent}

\def\build#1_#2^#3{\mathrel{\mathop{\kern 0pt#1}\limits_{#2}^{#3}}}

\def\date {le\ {\the\day}\ \ifcase\month\or janvier
\or fevrier\or mars\or avril\or mai\or juin\or juillet\or
ao\^ut\or septembre\or octobre\or novembre
\or d\'ecembre\fi\ {\oldstyle\the\year}}

\font\fivegoth=eufm5 \font\sevengoth=eufm7 \font\tengoth=eufm10

\newfam\gothfam \scriptscriptfont\gothfam=\fivegoth
\textfont\gothfam=\tengoth \scriptfont\gothfam=\sevengoth

\def\pro{\noindent {\bf Proof : }}

\def\smallsquare{\vbox{\hrule\hbox{\vrule height 1 ex\kern 1 ex\vrule}\hrule}}
\def\cqfd{\hfill \smallsquare\vskip 3mm}


\centerline{}

\vskip 4mm

\centerline{
\bf Palindromic continued fractions}

\vskip 8mm
\centerline{Boris A{\sevenrm DAMCZEWSKI} \ (Lyon) \
\& \ Yann B{\sevenrm UGEAUD} 
\footnote{*}{Supported by the Austrian Science
Fundation FWF, grant M822-N12. } 
\ (Strasbourg)}

\vskip 6mm

\vskip 8mm

\centerline{\bf 1. Introduction}

\vskip 6mm

It is widely believed that the continued fraction expansion of every
irrational algebraic number $\alpha$ 
is either eventually periodic (and this is
the case if, and only if, $\alpha$ is a quadratic irrationality) or
contains arbitrarily large partial quotients, but we seem to be
very far away from a proof (or a disproof). A preliminary step consists
in providing explicit examples of transcendental continued fractions.
The first result of this type is due to
Liouville \cite{Liouville}, who constructed real numbers
whose sequence of partial quotients grows very fast, too
fast for being algebraic. Subsequently, various authors
used deeper transcendence criteria from Diophantine
approximation to construct other classes of
transcendental continued fractions.
Of particular interest is the work of
Maillet \cite{Mai} (see also Section 34 of Perron 
\cite{Per}), who was the first to give examples of
transcendental continued fractions with bounded partial
quotients. 
Further examples have been provided by Baker \cite{Bak1,Bak2},
Davison \cite{Dav89}, Queff\'elec \cite{Queffelec98}, Allouche {\it et
al.} \cite{ADQZ} and
Adamczewski and Bugeaud \cite{Adamczewski_Bugeaud2},
among others.

\medskip

A common feature of the above quoted results 
is that they apply to real numbers whose continued fraction
expansion is `quasi-periodic', in the sense that 
it contains arbitrarily long blocks
of partial quotients which occur precociously at least twice.
In the present work, we investigate real numbers
whose sequence of partial quotients enjoys 
another combinatorial property, namely is `quasi-symmetrical', in the
sense that it begins in arbitrarily large `quasi-palindromes'. We provide
three new transendence criteria, that apply to a broad class of
continued fraction expansions, including expansions with unbounded
partial quotients. As in \cite{Adamczewski_Bugeaud2}, their 
proofs heavily depend on the Schmidt Subspace Theorem.
These criteria are stated in Section 2 and 
proved in Sections 6 and 7. Our method allows us to improve
upon an earlier result of A. Baker \cite{Bak1} on the 
transcendence of continued fractions whose sequence of
partial quotients is composed of long strings on $1$ and 
long strings of $2$, alternatively. See Section 3 for a
precise statement and Section 8 for its proof.
In Section 4, we provide an application of one of
our transcendence criteria to the explicit construction
of transcendental numbers with sharp properties of approximation by
rational numbers. All our auxiliary statements are gathered
in Section 5.

\vskip 18mm

\centerline{\bf 2. Main results}

\vskip 8mm

Throughout the present work, ${\cal A}$ 
denotes a given set, not necessarily finite. We identify any sequence
${\bf a}=(a_n)_{n \ge 1}$ of elements from ${\cal A}$
with the infinite word $a_1 a_2 \ldots a_n\ldots$
Recall that a finite word $a_1 a_2 \ldots a_n$ on ${\cal A}$ is
called a {\it palindrome} if $a_j = a_{n+1-j}$ for $j=1, \ldots, n$.

Our first transcendence criterion can be stated as follows.

\proclaim Theorem 1.
Let ${\bf a}=(a_n)_{n\geq 1}$ be a sequence of positive integers. 
If the word
${\bf a}$ begins in arbitrarily long palindromes, then the real number 
$\alpha:=[0;a_1, a_2, \ldots,a_n,\ldots]$ is either quadratic irrational or 
transcendental.

We point out that there is no assumption on the growth of the
sequence $(a_n)_{n\geq 1}$ in Theorem 1, unlike in Theorems 2 and 3 below.

As shown in \cite{AdBuMo}, 
given two distinct positive integers $a$ and $b$, 
Theorem 1 easily implies the transcendence of the
real number $[0; a_1, a_2, \ldots ...]$,
whose sequence of partial quotients is the Thue--Morse
sequence on the alphabet $\{a, b\}$, that is, with
$a_n = a$ (resp. $a_n = b$)
if the sum of binary digits of $n$ is odd (resp. even).
This result is originally due to M. Queff\'elec \cite{Queffelec98}.

\medskip

Before stating our next theorems, we need to introduce some more notation.
The length of a finite word
$W$ on the alphabet ${\cal A}$, that is, the number of letters
composing $W$, is denoted by $\vert W\vert$.
The mirror image of $W:= a_1 \ldots a_n$ is the word
$\uW := a_n \ldots a_1$. In particular, $W$ is a palindrome if,
and only if, $W = \uW$. A palindrome is thus a finite word invariant
under mirror symmetry. In order to relax this property of symmetry, we 
introduce the notion of {\it quasi-palindrome}. Let $U$ and $V$ be two
finite words; then, the word $UV\uU$ is called a quasi-palindrome of order
$w$, where $w=\vert V\vert/\vert U\vert$. Clearly, the larger
$w$ is, the weaker is the property of symmetry. In our next 
transcendence criterion, we replace the occurrences of aritrarily large
palindromes by the ones of arbitrarily large quasi-palindromes of 
bounded order. However, this weakening of our assumption has 
a cost, namely, an extra assumption 
on the growth of the partial quotients is then needed. Fortunately,
the latter assumption is not very restrictive. In particular, it is always 
satisfied by real numbers with bounded partial quotients. 

Let ${\bf a}=(a_n)_{n \ge 1}$ be a sequence of elements from ${\cal A}$.
Let $w$ be a rational number with $w>1$.
We say that ${\bf a}$ 
satisfies Condition $(*)_w$ if ${\bf a}$ is not
eventually periodic and if there exist 
two sequences of finite words 
$(U_n)_{n \ge 1}$ and $(V_n)_{n \ge 1}$ such that:

\medskip

\item{\rm (i)} For any $n \ge 1$, the word $U_n V_n \uU_n$ is a prefix
of the word ${\bf a}$;

\smallskip

\item{\rm (ii)} The sequence
$({\vert V_n\vert} / {\vert U_n\vert})_{n \ge 1}$ is bounded 
from above by $w$;

\smallskip

\item{\rm (iii)} The sequence $(\vert U_n\vert)_{n \ge 1}$ is 
increasing.

\medskip

We complement Theorem 1 in the following way.

\proclaim Theorem 2.
Let ${\bf a}=(a_n)_{n \ge 1}$ be a sequence of positive integers.
Let $(p_n/q_n)_{n \ge 1}$ denote the sequence of convergents to 
the real number
$$
\alpha:= [0; a_1, a_2, \ldots, a_n,\ldots].
$$ 
Assume that the sequence $(q_\ell^{1/\ell})_{\ell \ge 1}$ is bounded,
which is in particular the case when the sequence ${\bf a}$ is bounded.
If there exists a positive rational number $w$ such that
${\bf a}$ satisfies Condition $(*)_w$, 
then $\alpha$ is transcendental.

In the statements of Theorems 1 and 2 
the palindromes or the quasi-palindromes must appear at the very beginning 
of the continued fraction under consideration. Fortunately,
the ideas used in their proofs allow us to deal also 
with the more general situation where arbitrarily long quasi-palindromes 
occur not too far from the beginning.

Let $w$ and $w'$ be positive rational numbers.
We say that ${\bf a}$ 
satisfies Condition $(*)_{w,w'}$ if ${\bf a}$ is not
eventually periodic and if there exist 
three sequences of finite words $(U_n)_{n \ge 1}$,   
$(V_n)_{n \ge 1}$ and $(W_n)_{n \ge 1}$ such that:

\medskip

\item{\rm (i)} For any $n \ge 1$, the word $W_n U_n V_n \uU_n$ is a prefix
of the word ${\bf a}$;

\smallskip

\item{\rm (ii)} The sequence
$({\vert V_n\vert} / {\vert U_n\vert})_{n \ge 1}$ is bounded 
from above by $w$;

\smallskip

\item{\rm (iii)} The sequence
$({\vert U_n\vert} / {\vert W_n\vert})_{n \ge 1}$ is bounded 
from below by $w'$;

\smallskip

\item{\rm (iv)} The sequence $(\vert U_n\vert)_{n \ge 1}$ is 
increasing.

\medskip

We are now ready to complement Theorems 2 and 3 as follows.

\proclaim Theorem 3.
Let ${\bf a}=(a_n)_{n \ge 1}$ be a sequence of positive integers.
Let $(p_n/q_n)_{n \ge 1}$ denote the sequence of convergents to 
the real number
$$
\alpha:= [0; a_1, a_2, \ldots,a_n,\ldots].
$$ 
Assume that the sequence $(q_\ell^{1/\ell})_{\ell \ge 1}$ is bounded
and set $M = \limsup_{\ell \to + \infty} \, q_{\ell}^{1/\ell}$ and
$m = \liminf_{\ell \to + \infty} \, q_{\ell}^{1/\ell}$.
Let $w$ and $w'$ be positive rational numbers with
$$
w' >  2{\log M \over \log m} - 1. \eqno (2.1)
$$
If ${\bf a}$ satisfies Condition $(*)_{w,w'}$, 
then $\alpha$ is transcendental.

We display an immediate consequence of Theorem 3.

\proclaim Corollary 1.
Let ${\bf a}=(a_n)_{n \ge 1}$ be a sequence of positive integers.
Let $(p_n/q_n)_{n \ge 1}$ denote the sequence of convergents to 
the real number
$$
\alpha:= [0; a_1, a_2, \ldots,a_n,\ldots].
$$ 
Assume that the sequence $(q_{\ell}^{1/\ell})_{\ell \ge 1}$ converges.
Let $w$ and $w'$ be positive rational numbers with $w' > 1$. 
If ${\bf a}$ satisfies Condition $(*)_{w,w'}$, 
then $\alpha$ is transcendental.

Theorems 1 to 3 provide the exact analogues of Theorems 1 and 2
from \cite{Adamczewski_Bugeaud2}. It would be desirable to replace the
assumption (2.1) by the weaker one $w' > 0$. 
The statements of Theorems 2 and 3 show
that weakening the combinatorial assumption of Theorem 1
needs further assumptions on the size of the partial quotients.

\vskip 10mm

\centerline{\bf 3. On a theorem of A. Baker}

\vskip 6mm

In this Section we focus on a particular family of continued fractions 
introduced by Baker \cite{Bak1}. 
Let $a$ and $b$ be distinct positive integers.
Let ${\bf\Lambda}=(\lambda_n)_{n\geq 1}$ 
be a sequence of positive integers and set
$$
\alpha_{\bf\Lambda} :=[0;{\overline a}^{\lambda_1}, 
{\overline b}^{\lambda_2},{\overline a}^{\lambda_3},
\ldots,{\overline b}^{\lambda_{2n}}, 
{\overline a}^{\lambda_{2n+1}},\ldots],
$$
where, for positive integers $m$ and $\lambda$, we use the notation
${\overline m}^{\lambda}$ to denote a string of $\lambda$
consecutive partial quotients equal to $m$.
In his paper, Baker proved the transcendence of $\alpha_{\bf\Lambda}$ 
assuming that the sequence ${\bf\Lambda}$ increases
sufficiently rapidly.

For convenience, we assume that $b>a$ (this makes no restriction,
since the transcendence of a number does not depend on his first
partial quotients), and we set
$$
\alpha_a= [a; a, a, \ldots], \qquad \alpha_b = [b; b, b, \ldots],
$$
and
$$
\rho = \log \alpha_b / \log \alpha_a.
$$
Observe that we have $\rho > 1$.

Although stated in a weaker form, the following result
can be derived from \cite{Bak1}.

\proclaim Theorem. {\bf (A. Baker)}.
Let $\alpha_{\bf\Lambda}$ be as above. If the sequence ${\bf \Lambda}$ satisfies
$$
\liminf_{n \to \infty} \, {\lambda_{n+1} \over \lambda_n} >
{1 + \sqrt{8 \rho^2 + 1} \over 2 \rho},
$$
then $\alpha_{\bf\Lambda}$ is transcendental. 

Baker's proof rests on the generalisation, due to
LeVeque \cite{LeV}, of Roth's Theorem to
approximation by algebraic numbers from a given number field.
Using our approach based on quasi-palindromes and on the Schmidt
Subspace Theorem,
we are able to slightly improve upon his result.

\proclaim Theorem 4. Let $\alpha_{\bf\Lambda}$ be as above. 
If the sequence ${\bf \Lambda}$ satisfies
$$
\liminf_{n \to \infty} \, {\lambda_{n+1} \over \lambda_n} > \sqrt{2},
$$
then $\alpha_{\bf\Lambda}$ is transcendental.

Theorem 4 improves upon Baker's result, since we have
$$
\sqrt{2}< {1 + \sqrt{8 \rho^2 + 1} \over 2 \rho}<2,$$
for any $\rho > 1$. 
For instance, taking $a=1$ and $b=2$, Baker obtained the value 
$1.71 \ldots$ (cf. page 8 of \cite{Bak1}). 
Throughout the proof of Theorem 4 (postponed to Section 8), we will see
how it is often possible to refine the 
assumption (2.1) of Theorem 3.


\vskip 8mm

\centerline{\bf 4. Transcendental numbers with prescribed order of 
approximation}

\vskip 5mm

In Satz 6 of \cite{Ja31}, Jarn\'\i k
used the continued fraction theory to prove the
existence of real numbers with prescribed order of approximation by
rational numbers. Let 
$\varphi : \R_{\ge 1} \to \R_{>0}$ be a positive function.
We say that a real number $\alpha$
is `approximable at order $\varphi$' if there
exist infinitely many rational numbers $p/q$ with $q>0$ and
$|\alpha - p/q| < \varphi(q)$. 
Jarn\'\i k's result can then be stated as follows.

\proclaim Theorem J. Let $\varphi : \R_{\ge 1} \to \R_{>0}$ be a
non-increasing function such that $\varphi(x) = o(x^{-2})$ as $x$ tends
to infinity. Then, there are real numbers $\alpha$ which are
approximable at order $\varphi$ but which are not approximable at any order
$c \, \varphi$, with $0 < c < 1$.

In his proof, Jarn\'\i k constructed inductively the sequence of
partial quotients of $\alpha$. Actually, he showed that there are uncountably
many real numbers $\alpha$ with the required property, thus, in particular, 
transcendental numbers. However, his
construction did not provide any {\it explicit} example of such 
a transcendental $\alpha$.

In the present Section, we apply our Theorem 1 to
get, under an extra assumption on the function $\varphi$, explicit
examples of transcendental numbers satisfying the conclusion
of Theorem J. 

\proclaim Theorem 5. Let $\varphi : \R_{\ge 1} \to \R_{>0}$ be 
such that $x \mapsto x^2 \, \varphi(x)$ 
is non-increasing and tends to $0$ as $x$ tends
to infinity. Then, we can construct explicit examples
of transcendental numbers $\alpha$ which are
approximable at order $\varphi$ but which are not approximable at any order
$c \, \varphi$, with $0 < c < 1$.

\pro
Throughout the proof, for any real number $x$,
we denote by $\lceil x \rceil$ the smallest integer greater than
or equal to $x$.
We will construct inductively the sequence 
$(b_n)_{n \ge 1}$ of partial quotients of a suitable real number $\alpha$.
Denoting by $(p_n/q_n)_{n \ge 0}$ the sequence of convergents to
$\alpha$, it follows from the continued fraction theory
that, for any $n\ge 1$, we have
$$
{1 \over q_{n-1}^2 (b_n + 2)} < \biggl| \alpha - {p_{n-1} \over q_{n-1}}
\biggr| < {1 \over q_{n-1}^2 b_n}.  \eqno (4.1)
$$
Recall that $q_n \ge (3/2)^n$ for any $n\ge 5$. 
For any $x\ge 1$, set $\Psi (x) = x^2 \, \varphi(x)$.
Let $n_1 \ge 6$ be such that $\Psi((3/2)^{n}) \le 10^{-1}$ 
for any $n \ge n_1-1$. Then, set $b_1 = \ldots = b_{n_1 - 1} = 1$ and 
$b_{n_1} = \lceil 1/\Psi(q_{n_1-1}) \rceil$. Observe that $b_{n_1} \ge 10$.
Let $n_2 > n_1$ be such that $\Psi((3/2)^{n}) \le (10 b_{n_1})^{-1}$ 
for any $n \ge n_2-1$.
Then, set $b_{n_1 + 1} = \ldots = b_{n_2 - 1} = 1$ and 
$b_{n_2} = \lceil 1/\Psi(q_{n_2-1}) \rceil$. 
Observe that $b_{n_2} \ge 10 b_{n_1}$.

At this step, we have
$$
\alpha = [0; {\overline 1}^{n_1 - 1}, b_{n_1}, 
{\overline 1}^{n_2 - n_1-1}, b_{n_2}, \ldots],
$$
where, as in the previous Section, we denote by ${\overline 1}^m$ a sequence
of $m$ consecutive partial quotients equal to $1$.
Then, we complete by symmetry, in such a way that the continued fraction
expansion of $\alpha$ begins with a palindrome:
$$
\alpha = [0; {\overline 1}^{n_1 - 1}, b_{n_1}, 
{\overline 1}^{n_2 - n_1-1}, b_{n_2}, 
{\overline 1}^{n_2 - n_1-1}, b_{n_1}, 
{\overline 1}^{n_1 - 1}, \ldots].
$$
At this stage, we have constructed the first $2 n_2 - 1$ partial quotients
of $\alpha$. 
Let $n_3 > 2 n_2$ be such that $\Psi((3/2)^{n}) \le (10 b_{n_2})^{-1}$ 
for any $n \ge n_3-1$.
Then, set $b_{2 n_2} = \ldots = b_{n_3 - 1} = 1$ and 
$b_{n_3} = \lceil 1/\Psi(q_{n_3-1}) \rceil$. 
Observe that $b_{n_3} \ge 10 b_{n_2}$.
Then, we again complete by symmetry, and we repeat our process in
order to define $n_4$, $b_{n_4}$, and so on.

Clearly, the real number constructed in this way begins with infinitely
many palindromes, thus it is either quadratic or 
transcendental by Theorem 1. Moreover, the assumption on the function 
$\varphi$ implies that $\alpha$ has unbounded partial quotients. It thus 
follows that it is transcendental. 
It remains for us to prove that it has the
required property of approximation.

By (4.1), for any $j \ge 1$, we have 
$$
{\varphi (q_{n_j-1}) \over 1 + 3 \, q_{n_j - 1}^2 \, \varphi (q_{n_j-1})} <
\biggl| \alpha - {p_{n_j - 1} \over q_{n_j - 1}} \biggr| < \varphi (q_{n_j-1}).
\eqno (4.2)
$$
Let $p_n/q_n$ with $n \ge n_2$
be a convergent to $\alpha$ not in the subsequence
$(p_{n_j - 1} / q_{n_j - 1})_{j \ge 1}$, and let $k$ be the integer
defined by $n_k - 1 < n < n_{k+1} - 1$. 
Then, by (4.1) and the way we defined the $b_{n_j}$, we have
$$
\eqalign{
\biggl| \alpha - {p_n \over q_n} \biggr| & > {1 \over q_n^2 (b_{n+1} + 2)}
\ge  {1 \over q_n^2 (b_{n_{k-1}} + 2)} \cr
& \ge  {1 \over 3 q_n^2 \, b_{n_{k-1}}} \ge 
{\varphi (q_n) \over 3 q_{n_k - 1}^2 \, 
\varphi(q_{n_k-1}) \, b_{n_{k-1}}}, \cr}  \eqno (4.3)
$$
since $x \mapsto x^2 \varphi(x)$ is non-increasing.
We then infer from (4.3) and
$$
b_{n_k - 1} \le { b_{n_k} \over 10} \le {11 \over 100} \cdot
{1 \over q_{n_k - 1}^2 \, \varphi(q_{n_k-1})}
$$
that
$$
\biggl| \alpha - {p_n \over q_n} \biggr| \ge 3 \, \varphi(q_n). \eqno (4.4)
$$

To conclude, we observe that
it follows from (4.2) that $\alpha$ is approximable at order $\varphi$, and
from (4.2), (4.4) and the fact that $\varphi$ is non-increasing that
$\alpha$ is not approximable at any order $c \varphi$ with $0 < c < 1$.
The proof of Theorem 5 is complete. \cqfd

\vskip 6mm

\centerline{\bf 5. Auxiliary results}

\vskip 6mm

The proofs of Theorems 2 to 4 depend on a deep result
from Diophantine approximation, namely the powerful
Schmidt Subspace Theorem, stated as Theorem B below.
However, we do not need the full force of this theorem
to prove our Theorem 1: the transcendence
criterion given by Theorem A is sufficient for our purpose.

\proclaim Theorem A. {\bf (W. M. Schmidt).} 
Let $\alpha$ be a real number, which is
neither rational, nor quadratic. If there exist a real number $w>3/2$
and infinitely many triples of integers $(p, q, r)$ with $q > 0$ such that
$$
\max\biggl\{\biggl|\alpha - {p \over q} \biggr|, 
\biggl|\alpha^2 - {r \over q} \biggr|\biggr\} < {1 \over q^{w}},
$$
then $\alpha$ is transcendental.

\pro See \cite{Schm67}. \cqfd

\proclaim Theorem B. ({\bf W. M. Schmidt).} 
Let $m \ge 2$ be an integer.
Let $L_1, \ldots, L_m$ be linearly independent linear forms in
${\bf x} = (x_1, \ldots, x_m)$ with algebraic coefficients.
Let $\eps$ be a positive real number.
Then, the set of solutions ${\bf x} = (x_1, \ldots, x_m)$ 
in $\Z^m$ to the inequality
$$
\vert L_1 ({\bf x}) \ldots L_m ({\bf x}) \vert  \le
(\max\{|x_1|, \ldots , |x_m|\})^{-\eps}
$$
lies in finitely many proper subspaces of $\Q^m$.

\pro See e.g. \cite{Schmidt72a} or \cite{Schmidt80}. \cqfd

For the reader convenience, we 
further recall some well-known results from
the theory of continued fractions, 
whose proofs can be found e.g. in the book of Perron \cite{Per}.
The seemingly innocent Lemma 1 appears to be crucial in the
proofs of Theorems 2 to 4.

\proclaim Lemma 1.
Let $\alpha = [0; a_1, a_2, \ldots]$ be a real number with convergents
$(p_\ell/q_\ell)_{\ell \ge 1}$. Then, for any $\ell \ge 2$, we have
$$
{q_{\ell - 1} \over q_{\ell} } = [0 ; a_\ell , a_{\ell-1}, \ldots , a_1].
$$

\proclaim Lemma 2.
Let $\alpha = [0; a_1, a_2, \ldots]$ and $\beta =
[0; b_1, b_2, \ldots]$ be real numbers. Let $n \ge 1$ such that
$a_i = b_i$ for any $i=1, \ldots, n$. We then have
$|\alpha - \beta| \le q_n^{-2}$, where $q_n$ denotes the denominator
of the $n$-th convergent to $\alpha$.

For positive integers $a_1, \ldots, a_m$, we denote 
by $K_m (a_1, \ldots, a_m)$ the denominator of the rational number
$[0; a_1, \ldots, a_m]$. It is commonly called a {\it continuant}.

\proclaim Lemma 3.
For any positive integers $a_1, \ldots, a_m$ and any integer $k$ with
$1 \le k \le m-1$, we have
$$
K_m (a_1, \ldots , a_m) = K_m (a_m, \ldots, a_1)
$$
and
$$
\eqalign{
K_k (a_1, \ldots, a_k) \cdot K_{m-k} (a_{k+1}, \ldots, a_m)
& \le K_m (a_1, \ldots , a_m) \cr
& \le 2 \, K_k (a_1, \ldots, a_k) \cdot K_{m-k} (a_{k+1}, \ldots, a_m). \cr}
$$

The following three results are needed for the proof of Theorem 4.

\proclaim Lemma 4. 
Let $\alpha:=[0;a_1,a_2,\ldots,a_n]=p_n/q_n$ 
be a rational number $\alpha$. Then, we have
$$
\log q_n =\sum_{m=1}^{n}\log([a_m;a_{m-1},\ldots,a_1]).
$$

\medskip

\noi{\bf Proof.} 
It follows from Lemma 1 that 
$$
\log q_n=
\sum_{m=1}^n\log\left({q_m\over q_{m-1}}\right)=
\sum_{m=1}^{n}\log([a_m;a_{m-1},\ldots,a_1]),
$$
where 
$q_0,q_1, \ldots,q_n$ denote the denominators of the convergents of 
$\alpha$.  \cqfd

\medskip

Throughout the rest of this Section, $\theta$ denotes the
Golden Ratio $(1+\sqrt{5})/2$.

\medskip

\proclaim Lemma 5. Let $\alpha$ be a real number and denote by 
$(p_n/q_n)_{n\geq 1}$ the sequence of its convergents. Then, for 
every positive integer $k$, we have
$$
q_k\ge {\theta^{k+1}\over 2\sqrt{5}}.
$$

\medskip

\noi{\bf Proof.} By definition of the convergents 
we have $q_1=1$, $q_2=a_0+1$ and 
$q_n=a_nq_{n-1}+q_{n-2}$, 
where $(a_n)_{n\ge 0}$ denotes the sequence of partial quotients of $\alpha$. 
We thus have $q_1\ge 1$, $q_2\ge 1$ and 
$q_k\ge q_{k-1}+q_{k-2}$ for every integer $ k \ge 2$. 
For any integer $k > 0$, set $s_k := \theta^{k+1} / (2 \sqrt{5})$.
Since $s_1 < s_2 < 1$ and $s_{k+2} = s_{k+1} + s_k$, we get immediately
that $q_k \ge s_k$ for any positive integer $k$. \cqfd

\medskip

\proclaim Lemma 6. 
Let $r$ be a positive integer. Let $(x^{(k)})_{1\leq k\leq r}$ and 
 $(y^{(k)})_{1\leq k\leq r}$ be two finite sequences of real numbers lying 
in $[1,+\infty)$. 
Denote by $[x_{1,k};x_{2,k},\ldots]$ (resp.  
$[y_{1,k};y_{2,k},\ldots]$) the continued fraction expansion of $x^{(k)}$ 
(resp. 
of  $y^{(k)}$). If, for every $k$ satisfying $1\leq k\leq r$, we have 
$$
x_{j,k}=y_{j,k}, \qquad \hbox{for $j=1, \ldots, k$},
$$
then,  
$$
\sum_{k=1}^r\bigl\vert \log(x^{(k)})-\log(y^{(k)})\bigr\vert
<{20\over\theta^3}.
$$

\medskip

\noi{\bf Proof.} 
For every $k$ with $1\leq k\leq r$, denote by $q_k$ the denominator of the 
$k$-th convergent of $x^{(k)}$. By assumption, $q_k$ is also the 
denominator of the $k$-th convergent of $y^{(k)}$.  
It thus follows from Lemma 2 that
$$
\bigl\vert x^{(k)}-y^{(k)}\bigr\vert<{1\over q_k^2},
$$
which implies that  
$$ 
1-{1\over q_k^2}\leq {x^{(k)}\over y^{(k)}} \le 1+{1\over q_k^2}.
$$
We obtain 
$$
\bigl\vert\log(x^{(k)})-\log(y^{(k)})\bigr\vert
<{1\over q_k^2}
$$
and it follows from Lemma 5 
$$
\bigl\vert\log(x^{(k)})-\log(y^{(k)})\bigr\vert<{20\over\theta^{2k+2}}.
$$
We thus get that
$$
\sum_{k=1}^r\bigl\vert \log(x^{(k)})-\log(y^{(k)})\bigr\vert<
{20\over \theta^2}\,\sum_{k=1}^{+\infty}{1\over\theta^{2k}}
={20\over \theta^3},
$$
concluding the proof. \cqfd

\vskip 8mm

\centerline{\bf 6. Proof of Theorem 1}

\vskip 8mm

Let $n$ be a positive integer.
Denote by $p_n/q_n$ the $n$-th convergent to $\alpha$, that is, 
$p_n/q_n = [0; a_1, a_2, \ldots, a_n]$. By the theory of 
continued fraction, we have
$$
M_n := \pmatrix{ q_n & q_{n-1} \cr p_n & p_{n-1} \cr} =
\pmatrix{ a_1 & 1 \cr 1 & 0 \cr}
\pmatrix{ a_2 & 1 \cr 1 & 0 \cr} \ldots
\pmatrix{ a_n & 1 \cr 1 & 0 \cr}.
$$
Since such a decomposition is unique,
the matrix $M_n$ is symmetrical if, and only if,
the word $a_1 a_2 \ldots a_n$ is a palindrome.
Assume that this is case. Then, we have $p_n = q_{n-1}$.
Recalling that
$$
\biggl| \alpha - {p_n \over q_n} \biggr| < {1 \over q_n^2}
\quad {\rm and} \quad 
\biggl| \alpha - {p_{n-1} \over q_{n-1}} \biggr| < {1 \over q_{n-1}^2},
$$
we infer from $0 < \alpha < 1$, $a_1 = a_n$,
$|p_n q_{n-1} - p_{n-1} q_n| = 1$ and
$q_n \le (a_n + 1) q_{n-1}$ that
$$
\eqalign{
\biggl| \alpha^2 - {p_{n-1} \over q_n} \biggr| & \le
\biggl| \alpha^2 - {p_{n-1} \over q_{n-1}} \cdot  {p_n \over q_n} \biggr| 
\le \biggl| \alpha  + {p_{n-1} \over q_{n-1}} \biggr|
\cdot \biggl| \alpha -  {p_n \over q_n}  \biggr|
+ {1 \over q_n q_{n-1}} \cr
& \le 2 \biggl| \alpha  - {p_n \over q_n} \biggr| + {1 \over q_n q_{n-1}}
< {a_1 + 3 \over q_n^2}, \cr}
$$
whence
$$
\max\biggl\{\biggl|\alpha - {p_n \over q_n} \biggr|, 
\biggl|\alpha ^2 - {p_{n-1} \over q_n} \biggr|\biggr\} 
< {a_1 + 3 \over q_n^2}.  \eqno (6.1)
$$
Consequently, if the sequence of the partial quotients of $\alpha$ 
begins in arbitrarily long palindromes, then (6.1)
is satisfied for infinitely many integer triples
$(p_n, q_n, p_{n-1})$. By Theorem A, this shows that $\alpha$ is either
quadratic or transcendental.  \cqfd

\vskip 15mm

\centerline{\bf 7. Proofs of Theorems 2 and 3}

\vskip 6mm

Throughout the proofs of Theorems 2 and 3, for any finite word
$U = u_1 \ldots u_n$ on $\Z_{\ge 1}$, we denote by $[0; U]$ the
rational number $[0; u_1, \ldots , u_n]$.

\medskip

\noi{\bf Proof of Theorem 2.} 
Keep the notation and the hypothesis of this theorem.
Assume that the parameter $w > 1$ is fixed, as well as the 
sequences $(U_n)_{n \ge 1}$ and $(V_n)_{n \ge 1}$ occurring in the
definition of Condition $(*)_w$. 
Set also $r_n=\vert U_n \vert$ and $s_n=\vert U_n V_n \uU_n \vert$, 
for any $n \ge 1$.
We want to prove that the real number
$$
\alpha:= [0; a_1, a_2, \ldots]
$$ 
is transcendental. By assumption, we already know that 
$\alpha$ is irrational and not quadratic.
Therefore, we assume that $\alpha$ is algebraic of degree at
least three and we aim at deriving a contradiction.

Let $(p_{\ell}/q_{\ell})_{\ell \ge 1}$ denote the sequence of convergents 
to $\alpha$. 
The key fact for the proof of Theorem 2 is the equality
$$
{q_{\ell-1} \over q_\ell} = [0;a_\ell, a_{\ell-1}, \ldots,a_1],
$$ 
given by Lemma 1.
In other words, if $W_{\ell}$ denotes the prefix of 
length $\ell$ of the 
sequence ${\bf a}$, then
$q_{\ell-1}/q_\ell=[0; \overline{W_\ell}]$. 
Since, by assumption, we have 
$$
{p_{s_n} \over q_{s_n}}=[0;U_nV_n\overline{U_n}],
$$
we get that
$$
{q_{s_n-1} \over q_{s_n}}=[0;U_n\overline{V_n}\;\overline{U_n}],
$$ 
and it follows from Lemma 2 that 
$$
\vert q_{s_n}\alpha-q_{s_n-1}\vert<q_{s_n}q_{r_n}^{-2}.  \eqno (7.1)
$$
This shows in particular that
$$
\lim_{n \to + \infty} \, {q_{s_n-1} \over q_{s_n}} = \alpha. \eqno (7.2)
$$
Furthermore, we clearly have
$$
|q_{s_n} \alpha - p_{s_n}| < q_{s_n}^{-1}  \quad {\rm and} \quad
|q_{s_n-1} \alpha - p_{s_n-1}| < q_{s_n-1}^{-1}.
\eqno (7.3)
$$

Consider now the four linearly independent linear forms
with algebraic coefficients:
$$
\eqalign{
L_1(X_1, X_2, X_3, X_4) = & \alpha X_1 - X_3,  \cr
L_2(X_1, X_2, X_3, X_4) = & \alpha X_2 - X_4, \cr
L_3(X_1, X_2, X_3, X_4) = & \alpha X_1 - X_2, \cr
L_4(X_1, X_2, X_3, X_4) = & X_2. \cr}
$$
Evaluating them on the quadruple 
$(q_{s_n}, q_{s_n-1}, p_{s_n}, p_{s_n-1})$, it follows from (7.1)
and (7.3) that
$$
\prod_{1 \le j \le 4} \, |L_j (q_{s_n}, q_{s_n-1}, p_{s_n}, p_{s_n-1})|
< q_{r_n}^{-2}. \eqno (7.4)
$$
By assumption, there exists a real number $M$ such that
$$
\sqrt {2} \le q_{\ell}^{1/\ell} \le M
$$
for any integer $\ell \ge 3$. Thus, for any integer $n \ge 3$, we have
$$
q_{r_n}\ge \sqrt{2}^{r_n}
\ge (M^{s_n})^{(r_n \log  \sqrt{2})/(s_n \log M)} 
\ge q_{s_n}^{(r_n \log  \sqrt{2})/(s_n \log M)}\,
$$
and we infer from (7.4) and from $(ii)$ of Condition $(*)_w$ that
$$
\prod_{1 \le j \le 4} \, |L_j (q_{s_n}, q_{s_n-1}, p_{s_n}, p_{s_n-1})|
\ll q_{s_n}^{-\varepsilon}
$$
holds for some positive real number $\varepsilon$.

It then follows from
Theorem B that the points $(q_{s_n}, q_{s_n-1}, p_{s_n}, p_{s_n-1})$
lie in a finite number of proper subspaces of $\Q^4$. 
Thus, there exist a non-zero integer quadruple $(x_1,x_2,x_3,x_4)$ and
an infinite set of distinct positive integers ${\cal N}_1$ such that
$$
x_1 q_{s_n} + x_2 q_{s_n - 1} + x_3 p_{s_n} + x_4 p_{s_n - 1} = 0,  \eqno (7.5)
$$
for any $n$ in ${\cal N}_1$.
Dividing (7.5) by $q_{s_n}$, we obtain
$$
x_1 + x_2 {q_{s_n - 1} \over q_{s_n}} + x_3 {p_{s_n} \over q_{s_n}} 
+ x_4 { p_{s_n - 1} \over q_{s_n - 1} } \cdot {q_{s_n - 1} \over q_{s_n}}
= 0.  \eqno (7.6)
$$
By letting $n$ tend to infinity along ${\cal N}_1$ in (7.6),
it follows from (7.2) that 
$$
x_1 + (x_2 + x_3) \alpha + x_4 \alpha^2=0.
$$
Since, by assumption, $\alpha$ is not a quadratic number, 
we have $x_1=x_4=0$ and $x_2=-x_3$. Then, (7.5)
implies that 
$$
q_{s_n-1}=p_{s_n}. \eqno (7.7)
$$
Consider now the three linearly independent
linear forms with algebraic coefficients:
$$
L'_1(Y_1, Y_2, Y_3) = \alpha Y_1  - Y_2, \quad 
L'_2(Y_1, Y_2, Y_3) = \alpha Y_2 - Y_3, \quad 
L'_3(Y_1, Y_2, Y_3) = Y_1. 
$$
Evaluating them on the triple 
$(q_{s_n}, p_{s_n}, p_{s_n-1})$, we infer from (7.3) and
(7.7) that
$$
\prod_{1 \le j \le 3} \, |L'_j (q_{s_n}, p_{s_n}, p_{s_n-1})|
<  q_{s_n-1}^{-1} \ll q_{s_n}^{-0.9},
$$
since we have
$$
q_{\ell + 1} \ll q_{\ell}^{1.1}, \qquad
\hbox{for any $\ell \ge 1$,}  
$$
by Roth's Theorem. Here, the constants implied by $\ll$ depend
only on $\alpha$.

It then follows from
Theorem B that the points $(q_{s_n}, p_{s_n}, p_{s_n-1})$ 
with $n$ in ${\cal N}_1$ lie in a finite number of proper subspaces of 
$\Q^3$. 
Thus, there exist a non-zero integer triple $(y_1, y_2, y_3)$ and
an infinite set of distinct positive integers ${\cal N}_2$ such that
$$
y_1 q_{s_n} + y_2 p_{s_n} + y_3 p_{s_n-1}  = 0,  \eqno (7.8)
$$
for any $n$ in ${\cal N}_2$.
Dividing (7.8) by $q_{s_n}$, we get
$$
y_1 + y_2 {p_{s_n} \over q_{s_n}} + 
y_3 {p_{s_n-1} \over q_{s_n-1}} \cdot {q_{s_n-1} \over q_{s_n}} = 0.  
\eqno (7.9)
$$
By letting $n$ tend to infinity along ${\cal N}_2$, it thus follows 
from (7.7) and (7.9) that 
$$
y_1 + y_2 \alpha + y_3 \alpha^2=0.
$$
Since $(y_1,y_2,y_3)$ is a non-zero triple of integers, we
have reached a contradiction. Consequently, the real number $\alpha$ 
is transcendental. 
This completes the proof of the theorem. \cqfd

\bigskip

\noi{\bf Proof of Theorem 3.} 
Keep the notation and the hypothesis of this theorem.
Assume that the parameters $w$ and $w'$ are fixed, as well as the 
sequences $(U_n)_{n \ge 1}$, $(V_n)_{n \ge 1}$ and $(W_n)_{n \ge 1}$. 
Set also $r_n=\vert W_n\vert$, $s_n=\vert W_n U_n\vert$ and 
$t_n=\vert W_n U_n V_n \overline{U_n}\vert$, 
for any $n \ge 1$.
We want to prove that the real number
$$
\alpha:= [0; a_0, a_1, a_2, \ldots]
$$  
is transcendental.
By assumption, we already know that 
$\alpha$ is irrational and not quadratic.
Therefore, we assume that $\alpha$ is algebraic of degree at
least three and we aim at deriving a contradiction. 
Throughout this Section, the constants implied by $\ll$ depend
only on $\alpha$. In view of Theorem 2, we may assume that
$r_n \ge 1$ for any $n$.

The key idea for our proof is to consider, for any 
positive integer $n$, the 
rational $P_n/Q_n$ defined by
$$
{P_n \over Q_n}:=[0;W_n U_n V_n \overline{U_n}\,\overline{W_n}]
$$
and to use the fact that the word $W_n U_n V_n \overline{U_n}\,\overline{W_n}$ 
is a quasi-palindrome. 
Let $P'_n/Q'_n$ denote the last convergent to $P_n/Q_n$.     
By assumption we have 
$$
{p_{t_n} \over q_{t_n}} =[0;W_n U_n V_n \overline{U_n}]
$$
and it thus follows from Lemma 2 that 
$$
\vert Q_n\alpha-P_n\vert<Q_nq_{t_n}^{-2} \eqno (7.10)
$$
and 
$$
\vert Q'_n\alpha-P'_n\vert<Q'_n q_{t_n}^{-2}, \eqno (7.11)
$$
since $\uW_n$ has at least one letter.
Furthermore, Lemma 1 implies that
$$
{Q'_n \over Q_n}=[0;W_n U_n \overline{V_n}\,\overline{U_n}\,\overline{W_n}],
$$
and we get from  Lemma 2 that 
$$
\vert Q_n\alpha-Q'_n\vert<Q_n q_{s_n}^{-2}.  \eqno (7.12)
$$
This shows in particular that
$$
\lim_{n \to + \infty} \, {Q'_n \over Q_n} = \alpha. \eqno (7.13)
$$

Consider now the following 
four linearly independent linear forms with algebraic 
coefficients:
$$
\eqalign{
L_1(X_1, X_2, X_3, X_4) = & \alpha X_1 - X_3,  \cr
L_2(X_1, X_2, X_3, X_4) = & \alpha X_2 - X_4, \cr
L_3(X_1, X_2, X_3, X_4) = & \alpha X_1 - X_2, \cr
L_4(X_1, X_2, X_3, X_4) = & X_2. \cr}
$$
Evaluating them on the quadruple 
$(Q_n, Q'_n, P_n, P'_n)$, it follows from (7.10), (7.11) and (7.12) that
$$
\prod_{1 \le j \le 4} \, |L_j (Q_n, Q'_n, P_n, P'_n)|
< Q_n^4 q_{t_n}^{-4} q_{s_n}^{-2}.  \eqno (7.14)
$$

We infer from Lemma 3 that
$$
q_{t_n}q_{r_n} \le Q_n \le 2 q_{t_n}q_{r_n} \quad
{\rm and} \quad q_{s_n}^2 \le Q_n \le q_{t_n}^2, \eqno (7.15)
$$ 
and thus (7.14) gives  
$$
\prod_{1 \le j \le 4} \, |L_j (Q_n, Q'_n, P_n, P'_n)|
\ll q_{r_n}^4 \, q_{s_n}^{-2}. \eqno (7.16)
$$
Moreover, by our assumption (2.1),
there exists $\eta >0$ such that, for any $n$ large enough, we have 
$$
|U_n| \ge \biggl( 2 \, {\log M \over \log m} \cdot {1 + \eta \over 1 - \eta}
- 1 \biggr) \, |W_n|.
$$
Consequently, assuming that $n$ is sufficiently large, we get
$$
m^{(1- \eta) s_n} \ge M^{2(1+ \eta) r_n}
$$
and
$$
q_{s_n} \ge q_{r_n}^{2 + \eta'},
$$
for some positive real number $\eta'$.
It then follows from (7.16) that 
$$
\prod_{1 \le j \le 4} \, |L_j (Q_n, Q'_n, P_n, P'_n)|
\ll q_{s_n}^{-2 \eta' / (2 + \eta')}. \eqno (7.17)
$$
By assumption, we have for any $\ell$ large enough 
$$
\sqrt {2} \le q_{\ell}^{1/\ell} \le 2 M.
$$
Thus, for any integer $n$ large enough, we have
$$
\eqalign{
q_{s_n}\ge \sqrt{2}^{s_n}
\ge ((2M)^{t_n})^{(s_n \log  \sqrt{2})/( t_n \log 2M)} 
& \ge q_{t_n}^{(s_n \log  \sqrt{2})/(t_n \log 2M)} \cr
& \ge Q_n^{(s_n \log  \sqrt{2})/(2 t_n \log 2M)}, \cr}
$$
by (7.15).
We then infer from (7.17) and from $(ii)$ of Condition $(*)_{w, w'}$ that  
$$
\prod_{1 \le j \le 4} \, |L_j (Q_n, Q'_n, P_n, P'_n)|
\ll Q_n^{-\varepsilon}  \eqno (7.18)
$$
holds for some positive $\varepsilon$.

It then follows from
Theorem B that the points $(Q_n, Q'_n, P_n, P'_n)$
lie in a finite number of proper subspaces of $\Q^4$. 
Thus, there exist a non-zero integer quadruple $(x_1, x_2, x_3, x_4)$ and
an infinite set of distinct positive integers ${\cal N}_1$ such that
$$
x_1 Q_n + x_2 Q'_n + x_3 P_n + x_4 P'_n= 0,  \eqno (7.19)
$$
for any $n$ in ${\cal N}_1$. 
Dividing by $Q_n$, we obtain 
$$
x_1+x_2\,{Q'_n \over Q_n}+ x_3 \, {P_n \over Q_n} + 
x_4 \,{P'_n \over Q'_n} \cdot {Q'_n \over Q_n}= 0.
$$
By letting $n$ tend to infinity along ${\cal N}_1$,  
we infer from (7.13) that
$$
x_1 + (x_2 + x_3) \alpha + x_4 \alpha^2= 0.
$$ 
Since $(x_1, x_2, x_3, x_4)\not=(0,0,0,0)$ and since $\alpha$ is irrational
and not quadratic, we have $x_1 = x_4 = 0$ and $x_2 = -x_3$.
Then, (7.19) implies that 
$$
Q'_n=P_n. \eqno (7.20)
$$
Consider now the following 
three linearly independent linear forms with algebraic 
coefficients:
$$
L'_1(Y_1, Y_2, Y_3)=\alpha Y_1-Y_2,\;
L'_2(Y_1, Y_2, Y_3)=\alpha Y_2 - Y_3,\; 
L'_3(Y_1, Y_2, Y_3)=Y_1.
$$
Evaluating them on the quadruple 
$(Q_n, P_n, P'_n)$, it follows from (7.10), (7.11), (7.15) and (7.20) that
$$
\prod_{1 \le j \le 3} \, |L_j (Q_n, P_n, P'_n)|
\ll Q_n^3 q_{t_n}^{-4} \ll q_{r_n}^4 Q_n^{-1}
\ll q_{r_n}^4 q_{s_n}^{-2} \ll Q_n^{-\eps}, 
$$
with the same $\eps$ as in (7.18).
It then follows from
Theorem B that the points $(Q_n, P_n, P'_n)$
lie in a finite number of proper subspaces of $\Q^3$. 
Thus, there exist a non-zero integer triple $(y_1, y_2, y_3)$ and
an infinite set of distinct positive integers ${\cal N}_2$ such that
$$
y_1 Q_n + y_2 P_n + y_3 P'_n= 0,  
$$
for any $n$ in ${\cal N}_2$. 
We then proceed exactly as at the end of the proof of Theorem 2
to reach a contradiction. This finishes the proof of our theorem. \cqfd

\vskip 6mm

\centerline{\bf 8. Proof of Theorem 4}

\vskip 6mm

This Section is devoted to the proof of Theorem 4. 
Instead of a straightforward 
application of Theorem 3, which would yield a
weaker result, we carefully determine here  
the growth of the denominators of the convergents 
to the real numbers under consideration. 

\bigskip

\noi{\bf Proof of Theorem 4.} Let $\alpha= \alpha_{{\bf \Lambda}}=
[0;{\overline a}^{\lambda_1}, {\overline b}^{\lambda_2},
{\overline a}^{\lambda_3},\ldots ,{\overline a}^{\lambda_{2n-1}}
,{\overline b}^{\lambda_{2n}},\ldots]$ and denote by
$(p_\ell/q_\ell)_{\ell \geq 1}$ 
the sequence of its convergents. Assume that 
$$\liminf_{n\to\infty}{\lambda_{n+1}\over
\lambda_n}>\sqrt{2}. \eqno(8.1)$$ 

We first remark that if moreover the sequence 
$(\lambda_{n+1}/\lambda_n)_{n\geq
1}$ is not bounded from above, then the transcendence of $\alpha$
follows from a direct application of Theorem 3. 
Indeed, let us assume that there exists an increasing sequence of 
positive integers 
$(n_k)_{k\ge 1}$ such that 
$\lim_{k\to\infty}(\lambda_{n_k+1}/\lambda_{n_k})=+\infty$. Without loss 
of generality, we can also assume that $n_k$ is odd for every $k\ge 1$. Then, 
we apply Theorem 3 with 
$W_k={\overline a}^{\lambda_1}{\overline b}^{\lambda_2}
{\overline a}^{\lambda_3}\ldots {\overline a}^{\lambda_{n_k}}$,
$U_k={\overline b}^{\lfloor\lambda_{n_k+1}\rfloor/2}$ 
and $V_k$ equals to the empty word, since 
$\vert V_k\vert/\vert U_k\vert=0$ and 
$\lim_{k\to\infty}\vert U_k\vert/\vert W_k\vert=+\infty$.

From now on, we thus
assume that 
$$\limsup_{n\to\infty}{\lambda_{n+1}\over\lambda_n}<+\infty. \eqno(8.2)$$
Without loss of generality, we can
assume that $b>a$. For $n\geq 1$,  we set 
$W_n={\overline a}^{\lambda_1}{\overline b}^{\lambda_2}
{\overline a}^{\lambda_3}\ldots {\overline b}^{\lambda_{2n-2}}$, 
$U_n={\overline a}^{\lambda_{2n-1}}{\overline b}^{\lambda_{2n}}$ 
and $V_n={\overline a}^{\lambda_{2n+1}}
{\overline b}^{\lambda_{2n+2}-\lambda_{2n}}$. Set also 
$r_n=\vert W_n\vert$, $s_n=\vert W_nU_n\vert$ and $t_n=\vert
W_nU_nV_n{\overline U_n}\vert$. By assumption, we have 
$$p_{t_n}/q_{t_n}=[0;W_nU_nV_n{\overline U_n}].$$ 
Following the idea introduced in the proof of Theorem 3, we consider 
the rational $P_n/Q_n$ defined by
$$P_n/Q_n=[0;W_nU_nV_n{\overline U_n}\overline{W_n}].$$ 
Then, using the fact that the word 
$W_nU_nV_n{\overline U_n} \overline{W_n}$ is a quasi-palindrome, 
we can mimic the first
steps of the proof 
of Theorem 3 (see in particular (7.16)). We obtain that $\alpha$ is 
transcendental as soon as we can
ensure the existence of a positive real number $\varepsilon$ such 
that 
$$
q_{r_n}^4q_{s_n}^{-2}\ll Q_n^{-\varepsilon}.\eqno (8.3)
$$
Here and below, the numerical constant implicit in $\ll$ does not
depend on $n$.
Moreover, we can deduce from (8.1) and (8.2) that the sequence 
$(\vert U_n V_n\vert/\vert W_n\vert)_{n\geq 1}$ is bounded from above
which implies 
the existence of a positive real number $\eta$ such that 
$$q_{r_n}\ll Q_n^{\eta}.$$ It thus follows from (8.3) that $\alpha$ is
transcendental as soon as we can
ensure the existence of a positive real number $\delta$ such 
that 
$$
q_{r_n}^4q_{s_n}^{-2}\ll q_{r_n}^{-\delta}.\eqno (8.4)
$$ 
Let $(a_\ell)_{\ell \geq 1}$ 
denote the sequence of partial quotients of $\alpha$, and,
for every $m\geq 1$, set $x_m = [a_m; a_{m-1}, \ldots, a_1]$.
It follows from Lemma 4 that
$$
\log q_{s_n}=\sum_{m=1}^{s_n} \log x_m=
\log q_{r_n} + \sum_{m=r_n+1}^{s_n} \log x_m.
$$
In view of (8.4), it only
remains for us to prove that there exists a positive
real number $\delta$ such that 
$$\sum_{m=r_n+1}^{s_n} \log x_m\geq (1+\delta)\sum_{m=1}^{r_n} \log
x_m, \eqno (8.5)
$$
for $n$ large enough.

In the sequel of the proof, we show that (8.5) holds. 
For $ n \ge 1$, set
$c_n=\sum_{k=1}^{n}\lambda_k$, $d_n=\sum_{k=1}^{n}\lambda_{2k}$ 
and $e_n=\sum_{k=1}^{n}\lambda_{2k-1}$. 
To simplify the exposition, we put $c_0 = d_0 = e_0 = 0$ and,
for $0 \le j \le n-2$, we set
$$
A_{2j+1} = \sum_{m=c_{2j} + 1}^{c_{2j+1}} \log x_m, \quad 
B_{2j+2} = \sum_{m=c_{2j+1} + 1}^{c_{2j+2}} \log x_m, 
$$
and, for $0 \le j \le n-2$, we set
$$
A'_{2j+1} = \sum_{m=r_n+e_j+1}^{r_n+e_{j+1} } \log x_m, \quad 
B'_{2j+2} = \sum_{m=r_n+\lambda_{2n-1}+d_j+1}^{r_n+\lambda_{2n-1}+
d_{j+1} } \log x_m. $$
It follows from (1) that for $0 \le j \le n-2$ and for 
$r_n+e_j+1\le m\le r_n+e_{j+1}$, 
the first $m$ partial quotients of $x_m$ are all 
equals to $a$. We similary deduce from (1) that for $0 \le j \le n-2$ and for 
$r_n+\lambda_{2n-1}+d_j+1\le m\le r_n+\lambda_{2n-1}+
d_{j+1}$, the first $m$ partial quotients of $x_m$ are all equal to $b$. 
By Lemma 6, we thus have
$$
\eqalign{
|A_{2j+1} - \log \alpha_a| < {20\over \theta^3}, 
\quad & |A'_{2j+1} - \log \alpha_a|
< {20\over \theta^3} , \cr
|B_{2j+2} - \log \alpha_b| < {20\over \theta^3}, 
\quad & |B'_{2j+2} - \log \alpha_b|
< {20\over \theta^3}, \cr}\eqno (8.6)
$$
where, as in Section 5, we have set $\theta = (1 + \sqrt{5})/2$.
On the other hand, we have
$$
\sum_{m=1}^{r_n} \log x_m=\sum_{j=0}^{n-2}(A_{2j+1}+B_{2j+2})
\eqno (8.7)
$$
and
$$\sum_{m=r_n+1}^{s_n} \log x_m=\sum_{j=0}^{n-2}(A'_{2j+1}+B'_{2j+2})+
\sum_{m=r_n+e_{n-1}+1}^{r_n+\lambda_{2n-1}}\log x_m+
\sum_{m=s_n-\lambda_{2n}+d_{n-1}+1}^{s_n}\log x_m.$$
Then,  (8.6) and (8.7) imply that
$$\sum_{m=r_n+1}^{s_n} \log x_m-\sum_{m=1}^{r_n} \log x_m\geq 
\sum_{m=r_n+e_{n-1}+1}^{r_n+\lambda_{2n-1}}\log x_m+
\sum_{m=s_n-\lambda_{2n}+d_{n-1}+1}^{s_n}\log x_m-{80n\over\theta^3}.
\eqno (8.8)$$

We also deduce from (8.1) the existence of a positive real number  
$\omega$ such that 
$$
\lambda_{2n-1}>(1+\omega)e_{n-1}\quad {\rm and} \quad  \lambda_{2n}
>(1+\omega)d_{n-1}. 
$$
Now, from (8.8) and the fact that $x_m \ge (b+2)/(b+1)$ for any $m\ge 2$, 
it thus follows 
$$
\sum_{m=r_n+1}^{s_n} \log x_m-\sum_{m=1}^{r_n} \log x_m\geq 
w(e_{n-1}+d_{n-1})\log{b+2 \over b+1}-{80n\over\theta^3}=
\omega r_n\log{b+2 \over b+1}-{80n\over\theta^3}.
$$
We easily deduce from  (8.1) that $r_n \ge 2^n$ for $n$ large enough, 
which implies that  
$$
\sum_{m=r_n+1}^{s_n} \log x_m-\sum_{m=1}^{r_n} \log x_m\geq {\omega\over 2}
r_n\log{b+2 \over b+1},\eqno (8.9)
$$
for $n$ large enough. On the other hand, we have $x_m<b+1$ for every
positive integer $m$, implying  that
$$
\sum_{m=1}^{r_n} \log x_m<r_n\log(b+1). \eqno (8.10)
$$
Combining (8.9) and (8.10), we get
the existence of a positive real number $\delta$ such that 
$$
\sum_{m=r_n+1}^{s_n} \log x_m-\sum_{m=1}^{r_n} \log x_m>
\delta\left(
\sum_{m=1}^{r_n} \log x_m\right),
$$
for $n$ large enough. In view of (8.5), this concludes the proof. \cqfd

\vskip 12mm

\centerline{\bf References}

\vskip 7mm

\beginthebibliography{999}

\bibitem{Adamczewski_Bugeaud2}
B. Adamczewski \& Y. Bugeaud,
{\it On the complexity of algebraic numbers. II.
Continued fractions}, Acta Math. To appear.

\bibitem{AdBuMo}
B. Adamczewski \& Y. Bugeaud,
{\it A short proof of the transcendence of the
Thue--Morse continued fraction}. Preprint.

\bibitem{ADQZ}
J.-P. Allouche, J. L. Davison, M. Queff\'elec \& L. Q. Zamboni,
{\it Transcendence of Sturmian or morphic continued fractions},
J. Number Theory 91 (2001), 39--66.

\bibitem{Bak1} 
A. Baker,
{\it Continued fractions of transcendental numbers},
Mathematika 9 (1962), 1--8.

\bibitem{Bak2} 
A. Baker,
{\it On Mahler's classification of transcendental numbers},
{Acta Math.}
{111}
{(1964)},
{97--120}.

\bibitem{Dav89}
J. L. Davison,
{\it A class of transcendental numbers with bounded partial quotients}.
In R. A. Mollin, ed., Number Theory and Applications, pp. 365--371, Kluwer
Academic Publishers, 1989.

\bibitem{Ja31}
V. Jarn\'\i k, 
{\it \"Uber die simultanen Diophantische Approximationen},
Math. Z. 33 (1931), 505--543.

\bibitem{Khintchine}
A. Ya. Khintchine,
Continued fractions, Gosudarstv. Izdat. Tehn.-Theor. Lit. 
Moscow-Leningrad, 2nd edition, 1949 (in Russian).

\bibitem{Lang}
S. Lang,
Introduction to Diophantine Approximations, Sprin\-ger-Verlag (1995).

\bibitem{LeV}
W. J. LeVeque,
Topics in Number Theory, Vols. 1 \& 2, 
Addison--Wesley Publ. Co., Inc., Reading, MA.

\bibitem{Liouville}
J. Liouville, 
{\it Sur des classes tr\`es \'etendues de quantit\'es dont la valeur 
n'est ni alg\'ebri\-que, ni m\^eme r\'eductible \`a des irrationelles 
alg\'ebriques},
C. R. Acad. Sci. Paris 18 (1844), 883--885 and 993--995.

\bibitem{Mai}
E. Maillet,
Introduction \`a la th\'eorie des nombres transcendants et des propri\'et\'es
arithm\'etiques des fonctions, Gauthier-Villars, Paris, 1906.

\bibitem{Per}
O. Perron,
Die Lehre von den Ketterbr\"uchen.
Teubner, Leipzig, 1929.

\bibitem{Queffelec98}
M. Queff\'elec,
{\it Transcendance des fractions continues de Thue--Morse},
J. Number Theory 73 (1998), 201--211.

\bibitem{Schm67}
 W. M. Schmidt,
{\it On simultaneous approximations of two algebraic numbers by rationals},
Acta Math. 119 (1967), 27--50.

\bibitem{Schmidt72a}
 W. M. Schmidt,
{\it Norm form equations},
Ann. of Math. 96 (1972), 526--551.
  
\bibitem{Schmidt80}
W. M. Schmidt, 
{\it Diophantine approximation},
Lecture Notes in Mathematics 785, Springer, Berlin, 1980.

\endthebibliography


\bigskip\bigskip

\noindent Boris Adamczewski   \hfill{Yann Bugeaud}

\noindent   CNRS, Institut Camille Jordan  
\hfill{Universit\'e Louis Pasteur}

\noindent   Universit\'e Claude Bernard Lyon 1 
\hfill{U. F. R. de math\'ematiques}

\noindent   B\^at. Braconnier, 21 avenue Claude Bernard
 \hfill{7, rue Ren\'e Descartes}

\noindent   69622 VILLEURANNE Cedex (FRANCE)   
\hfill{67084 STRASBOURG Cedex (FRANCE)}

\vskip2mm
 
\noindent {\tt Boris.Adamczewski@math.univ-lyon1.fr}
\hfill{{\tt bugeaud@math.u-strasbg.fr}}

\vskip2mm
\hfill{Institut f\"ur Diskrete Mathematik und Geometrie}

\hfill{TU Wien}

\hfill{Wiedner Hauptstrasse 8--10}

\hfill{1040 WIEN (AUSTRIA)}

\vskip 2mm

\hfill{{\tt bugeaud@geometrie.tuwien.ac.at}}

\bye

\bye

%% file: typpo
%
\catcode`@=11
%
%
\def\bibn@me{R\'ef\'erences}
\def\bibliographym@rk{\centerline{{\sc\bibn@me}}
	\sectionmark\section{\ignorespaces}{\unskip\bibn@me}
	\bigbreak\bgroup
	\ifx\ninepoint\undefined\relax\else\ninepoint\fi}
%
%
%
\let\refsp@ce=\ 
\let\bibleftm@rk=[
\let\bibrightm@rk=]
%
%
%
\def\numero{n\raise.82ex\hbox{$\fam0\scriptscriptstyle o$}~\ignorespaces}
%
%
\newcount\equationc@unt
\newcount\bibc@unt
\newif\ifref@changes\ref@changesfalse
\newif\ifpageref@changes\ref@changesfalse
\newif\ifbib@changes\bib@changesfalse
\newif\ifref@undefined\ref@undefinedfalse
\newif\ifpageref@undefined\ref@undefinedfalse
\newif\ifbib@undefined\bib@undefinedfalse
\newwrite\@auxout
%
%
\def\eqnum{\global\advance\equationc@unt by 1%
\edef\lastref{\number\equationc@unt}%
\eqno{(\lastref)}}
%
%
%
%
%
%
\def\re@dreferences#1#2{{%
	\re@dreferenceslist{#1}#2,\undefined\@@}}
\def\re@dreferenceslist#1#2,#3\@@{\def\next{#2}%
	\expandafter\ifx\csname#1@@\meaning\next\endcsname\relax
	??\immediate\write16
	{Warning, #1-reference "\next" on page \the\pageno\space
	is undefined.}%
	\global\csname#1@undefinedtrue\endcsname
	\else\csname#1@@\meaning\next\endcsname\fi
	\ifx#3\undefined\relax
	\else,\refsp@ce\re@dreferenceslist{#1}#3\@@\fi}
%
%
%
\def\newlabel#1#2{{\def\next{#1}\newl@bel#2}}
\def\newl@bel#1#2{%
	\expandafter\xdef\csname ref@@\meaning\next\endcsname{#1}%
	\expandafter\xdef\csname pageref@@\meaning\next\endcsname{#2}}
\def\label#1{{%
	\toks0={#1}\message{ref(\lastref) \the\toks0,}%
	\ignorespaces\immediate\write\@auxout%
	{\noexpand\newlabel{\the\toks0}{{\lastref}{\the\pageno}}}%
	\def\next{#1}%
	\expandafter\ifx\csname ref@@\meaning\next\endcsname\lastref%
	\else\global\ref@changestrue\fi%
	\newlabel{#1}{{\lastref}{\the\pageno}}}}
\def\ref#1{\re@dreferences{ref}{#1}}
\def\pageref#1{\re@dreferences{pageref}{#1}}
%
%
\def\bibcite#1#2{{\def\next{#1}%
	\expandafter\xdef\csname bib@@\meaning\next\endcsname{#2}}}
\def\cite#1{\bibleftm@rk\re@dreferences{bib}{#1}\bibrightm@rk}
%
%
\def\beginthebibliography#1{\bibliographym@rk
	\setbox0\hbox{\bibleftm@rk#1\bibrightm@rk\enspace}
	\parindent=\wd0
	\global\bibc@unt=0
	\def\bibitem##1{\global\advance\bibc@unt by 1
		\edef\lastref{\number\bibc@unt}
		{\toks0={##1}
		\message{bib[\lastref] \the\toks0,}%
		\immediate\write\@auxout
		{\noexpand\bibcite{\the\toks0}{\lastref}}}
		\def\next{##1}%
		\expandafter\ifx
		\csname bib@@\meaning\next\endcsname\lastref
		\else\global\bib@changestrue\fi%
		\bibcite{##1}{\lastref}
		\medbreak
		\item{\hfill\bibleftm@rk\lastref\bibrightm@rk}%
		}
	}
\def\endthebibliography{\egroup\par}
%
%
%
\def\@closeaux{\closeout\@auxout
	\ifref@changes\immediate\write16
	{Warning, changes in references.}\fi
	\ifpageref@changes\immediate\write16
	{Warning, changes in page references.}\fi
	\ifbib@changes\immediate\write16
	{Warning, changes in bibliography.}\fi
	\ifref@undefined\immediate\write16
	{Warning, references undefined.}\fi
	\ifpageref@undefined\immediate\write16
	{Warning, page references undefined.}\fi
	\ifbib@undefined\immediate\write16
	{Warning, citations undefined.}\fi}
%
%
\immediate\openin\@auxout=\jobname.aux
\ifeof\@auxout \immediate\write16
  {Creating file \jobname.aux}
\immediate\closein\@auxout
\immediate\openout\@auxout=\jobname.aux
\immediate\write\@auxout {\relax}%
\immediate\closeout\@auxout
\else\immediate\closein\@auxout\fi
%
%
\input\jobname.aux
\immediate\openout\@auxout=\jobname.aux
%
%
\catcode`@=12